\documentclass{amsart}
\usepackage{amssymb}
\usepackage{amscd,enumerate,amsfonts,calc,amsmath,verbatim}
                    
\newcommand{\Ext}{\operatorname{Ext}}

\newcommand{\Tor}{\operatorname{Tor}}
\newcommand{\Hom}{\operatorname{Hom}}

\newcommand{\Ima}{\operatorname{Im}}
\newcommand{\Tr}{\operatorname{Tr}}
\newcommand{\Coker}{\operatorname{Coker}}
\newcommand{\rank}{\operatorname{rank}}

\newcommand{\col}{\colon}

\newcommand{\ges}{\geqslant}
\newcommand{\Ker}{\operatorname{Ker}}

\newtheorem*{Theorem}{Theorem}

\newcommand{\m}{{\mathfrak m}}
\newcommand{\fm}{{\mathfrak m}}
\newcommand{\fn}{{\mathfrak n}}

\newcommand{\bd}{\boldsymbol}

\theoremstyle{remark}

\swapnumbers

\theoremstyle{plain}
\newtheorem{theorem}{Theorem}[section]

\newtheorem{proposition}[theorem]{Proposition}
\newtheorem{lemma}[theorem]{Lemma}

\theoremstyle{definition}
\newtheorem{chunk}[theorem]{}

\newtheorem*{chunk*}{}

\theoremstyle{remark}

\newtheorem{remark}[theorem]{Remark}
\newtheorem{remarks}[theorem]{Remarks}

\newtheorem*{question}{Question}

\numberwithin{equation}{theorem}

\numberwithin{subchunk}{theorem}

\begin{document}\title[Independence of Total Reflexivity Conditions]
{Independence of the total reflexivity \\ conditions for modules} 
\author[D.~A.~Jorgensen]{David A.~Jorgensen } 
\address{Department of Mathematics, University of Texas at Arlington, 
Arlington, TX 76019}
\email{djorgens@math.uta.edu}
\author[L.~M.~\c Sega]{Liana M.~\c Sega}
\address{Department of Mathematics, Michigan State University, 
East Lansing, MI 48824}
\email{lsega@math.msu.edu}

\date{\today} 

\begin{abstract}
  We show that the conditions defining total reflexivity for modules
  are independent. In particular, we construct a
  commutative Noetherian local ring $R$ and a reflexive $R$-module $M$
  such that $\Ext^i_R(M,R)=0$ for all $i>0$, but $\Ext^i_R(M^*,R)\ne
  0$ for all $i>0$.
\end{abstract}
\maketitle

\section*{introduction}

Let $R$ be a commutative Noetherian ring. For any $R$-module $M$ we
set $M^*=\Hom_R(M,R)$. An $R$-module $M$ is said to be {\it reflexive}
if it is finite and the canonical map $M\to M^{**}$ is bijective.  A
finite $R$-module $M$ is said to be {\it totally reflexive} if it
satisfies the following conditions:
\begin{enumerate}[\quad (i)]
\item $M$ is reflexive
\item $\Ext^i_R(M,R)=0$ for all $i>0$
\item $\Ext^i_R(M^*,R)=0$ for all $i>0$. 
\end{enumerate}

This notion is due to Auslander and Bridger \cite{ABr}: the totally
reflexive modules are precisely the modules of {\it G-dimension} zero.
The G-dimension of a module is one of the best studied non-classical
homological dimensions, and is defined in terms of the length of
a resolution of the module by totally reflexive modules.

Given any homological dimension, a serious concern is whether its
defining conditions can be verified effectively.  For example, the 
projective dimension of a finite $R$-module $M$ is zero 
if and only if $\Ext^1_R(M,N)=0$ for all finite
$R$-modules $N$.  However, when $R$ is local with maximal ideal $\m$,
one only needs to check
vanishing for $N=R/\m$.  In the same spirit, it is natural to ask
whether the set of conditions defining total reflexivity is
overdetermined (cf. \cite[\S 2]{AM}) and in particular, whether total
reflexivity for a module can be established by verifying vanishing of
only finitely many Ext modules.

When $R$ is a local Gorenstein ring, (ii) implies the other two
 conditions above, and it is equivalent to $M$ being maximal
 Cohen-Macaulay. Recently, Yoshino \cite{Y} studied other situations when (ii)
 alone implies total reflexivity, and raised the question whether this
 is always the case.

In the present paper we give an example of a local Artinian ring $R$
which admits modules whose total reflexivity conditions are
independent, in that (ii) implies neither (i) nor (iii); (i) and (ii)
do not imply (iii), equivalently, (i) and (iii) do not imply (ii).
More precisely, we prove the following result as Theorem \ref{main}:
\begin{Theorem}
\label{A}
There exists a local Artinian ring $R$, and a family 
$\{M_s\}_{s\ges 1}$ of reflexive $R$-modules such that
\begin{enumerate}[\quad\rm(1)]
\item $\Ext^i_R(M_s,R)=0$ if and only if $1\le i\leq s-1$;
\item $\Ext^i_R(M_s^*,R)=0$ for all $i>0$.
\end{enumerate}
Moreover, there exists a non-reflexive $R$-module $L$
such that
\begin{enumerate}[\quad\rm(1$'$)]
\item $\Ext^i_R(L,R)=0$ for all $i>0$;
\item $\Ext^i_R(L^*,R)\ne 0$ for all $i>0$. 
\end{enumerate}
\end{Theorem}

By taking $N_s=M^*_s$ for all $s\ge 1$,
we get a statement dual to that of the first
part above: there exits a family $\{N_s\}_{s\ges 1}$ of reflexive 
$R$-modules such that (1) $\Ext^i_R(N_s^*,R)=0$ if and only if 
$1\le i\leq s-1$;
(2) $\Ext^i_R(N_s,R)=0$ for all $i>0$.

This theorem shows that in order to check whether or not a module 
$M$ is totally reflexive --- even for a local Artinian ring ---
one needs to check vanishing of $\Ext^i_R(M,R)$ and
$\Ext^i_R(M^*,R)$ for infinitely many values of $i$.  
In Section 2, however, we point out that when $R$ is Artinian and standard 
graded, in the sense that $R=\bigoplus_{i=0}^{\infty}R_i$ with
$R_0=k$, a field, and $R=R_0[R_1]$, one may 
skip checking finitely many values of $i$ of the same parity.

In our example, $R$ is a standard graded Koszul algebra and has Hilbert series
$\sum_{i\ge 0}\rank_k(R_i)t^i=1+4t+3t^2$. The ring $R$ is thus local,
and its maximal ideal $\m$ satisfies $\m^3=0$. Note that our example
is minimal in the following sense: if $\fm^2=0$, then any finite
$R$-module $M$ which satisfies $\Ext^i_R(M,R)=0$ for some $i>1$ is
totally reflexive, hence (ii) alone implies total reflexivity. 
(See 1.1 below.)

Our construction involves a minimal acyclic complex C of finite free
$R$-modules such that the sequence $\{\rank_R(C_i)\}_{i\ges 0}$ is
strictly increasing and has exponential growth, while the sequence
$\{\rank_R(C_{-i})\}_{i\ges 0}$ is constant. In the last section we
raise several related questions. 

\section{Independence}

Let $R$ be a Noetherian commutative ring and $M$ a finite
$R$-module. Suppose that $\phi\col G\to F$ is a homomorphism of finite
free $R$-modules with $M=\Coker\phi$. The $R$-module $\Tr(M):=\Coker
\phi^*$ is called the {\it transpose} of $M$, and it is unique up to
projective equivalence; it is thus well-defined in the stable
category of $R$. The $R$-modules $\Tr(\Tr(M))$ and $M$ are
isomorphic up to projective summands and there is an exact sequence
$$0\to\Ext^1_R(\Tr(M),R)\to M\to M^{**}\to \Ext^2_R(\Tr(M),R)\to 0$$
where the map in the middle is the natural map. In particular, $M$ is
reflexive if and only if $\Ext^i_R(\Tr(M),R)=0$ for $i=1,2$. Also note
that $M^*$ is a second syzygy of $\Tr(M)$, hence $\Ext^i_R(M^*,R)\cong
\Ext^{i+2}_R(\Tr(M),R)$ for all $i>0$. Thus the definition of $M$
being totally reflexive can be recast as follows:
$$\Ext^i_R(M,R)=0=\Ext^i_R(\Tr(M),R)\quad\text{for all}\quad i>0\,.$$

It is convenient to have a uniform notation for the 
conditions above: let $i\in\mathbb Z$, $i\neq 0$. 
We say that $M$ satisfies condition $(\mathsf{TR}_i)$ provided 
\[
(\mathsf{TR}_i):\hskip .5in \begin{cases} \Ext^i_R(M,R)=0 &\text{ if $i\ge 1$} \\
\Ext^{-i}_R(\Tr(M),R)=0 &\text{ if $i\le -1$}.
\end{cases}
\]
Thus $M$ is totally reflexive if and only if $M$ satisfies
$(\mathsf{TR}_i)$ for all $i\ne 0$.

Note that we did not define a condition $(\mathsf{TR}_i)$ for
$i=0$. When referring to these conditions we will assume tacitly that
$i\ne 0$. 

\begin{remarks}
(1) Assume that $R$ is local Gorenstein. The module $M$ is then totally
    reflexive if and only if $(\mathsf{TR}_i)$ is satisfied for all
    $i$ with $0< i\le \dim R$.

(2)  Assume that $R$ is local, with maximal ideal $\fm$ and residue
     field $k$. If $\fm^2=0$, then the first syzygy $N$ in a minimal
     free resolution of $M$ has $\m N=0$, hence it is a finite
     dimensional $k$-vector space. If $\Ext^i_R(M,R)=0$ for some $i>1$
     then either $M$ is free or $\Ext^{i-1}_R(k,R)=0$. Thus, when
     $\fm^2=0$ we obtain: $M$ is totally reflexive if and only if
     $(\mathsf{TR}_i)$ is satisfied for some $i$ with $i>1$ or $i<-1$.

(3) Assume that $R$ is a commutative Noetherian ring.  Yoshino 
 proved in \cite{Y} that if the full subcategory of $R$-modules
    $M$ with
    $\Ext^i_R(M,R)=0$ for all $i>0$ is of finite type, then a module
    $M$ is totally reflexive if and only if it is an object of this 
    subcategory. 
\end{remarks}

The remarks above lead to the question: How many conditions
$(\mathsf{TR}_i)$ does one need to check for total reflexivity,
without placing extra assumptions on the ring?  Is it possible that
only finitely many suffice? Is it enough to check $(\mathsf{TR}_i)$
for all $i>0$, or, more generally, for all $i>s$ for some integer $s$?

Theorem \ref{main} (stated also in the
introduction) provides negative answers to these questions: in
general, one needs to check the conditions $(\mathsf{TR}_i)$ for
infinitely many positive values of $i$ and infinitely many negative
values of $i$.

We now describe the ring of Theorem \ref{main}. Related rings were
used in \cite{GP} and then in \cite{JS} to disprove various
conjectures.
\medskip

Let $k$ be a field which is not algebraic over a finite field and let
$\alpha\in k$ be an element of infinite multiplicative order.  For the
remainder of this section we assume the ring $R$ to be defined as
follows.

\begin{chunk} \label{ringR} Consider the polynomial ring
  $k[V,X,Y,Z]$ in four variables (each of degree one) and set 
\[
R=k[V,X,Y,Z]/I,
\]
  where $I$ is the ideal generated by the following quadratic
  relations:
\[
V^2,\,Z^2,\,XY,\,VX+\alpha XZ,\,VY+YZ,\,VX+Y^2,\,VY-X^2\,.
\]
As a vector space over $k$, it has a basis consisting of the following
$8$ elements:
\[
1,\,v,\,x,\,y,\,z,\,vx,\,vy,\,vz,
\]
where $v,x,y,z$ denote the residue classes of the variables modulo
$I$. In particular, $R$ has Hilbert series $H_R(t)=1+4t+3t^2$.

\begin{remark}
One may check that the generators for $I$ listed above are a 
Gr\"obner basis for $I$.  Therefore
by \cite[Section 4]{F}, the ring $R$ is Koszul, and it follows that the
Poincar\'e series of its residue field $k$ is
$\sum_{i}\rank_k\Tor_i^R(k,k)t^i=(1-4t+3t^2)^{-1}$.
\end{remark}
\end{chunk}

\smallskip

For each integer $i\leq 0$ we let $d_i\colon R^2\to R^2$ denote the map
given with respect to the standard basis of $R^2$ by the matrix
\[
\left(\begin{matrix}
v& \alpha^{-i}x\\
y& z\end{matrix}\right).
\]
Also, let $d_1\col R^3\to R^2$ denote the map represented by the
matrix 
$$
\left(\begin{matrix}v& \alpha^{-1}x& yz\,\\
y& z&0\end{matrix}\right)
$$ 
and let $d_2\col R^7\to R^3$ be represented by the matrix
$$
\left(\begin{matrix}v& \alpha^{-2}x& -y&0&0&0&0\,\\
y& z&\alpha x&0&0&0&0\\0&0&0&v&x&y&z\end{matrix}\right).
$$ 
 
Consider a minimal free resolution of $\Coker d_2$ with $d_2$ as
the first differential:
$$
\cdots \to R^{b_i}\xrightarrow{d_i} R^{b_{i-1}}\to \cdots\to R^{b_3}\xrightarrow{d_3}R^7\xrightarrow{d_2}R^3,
$$ 
where for each $i\ge 3$ the map $d_i\col R^{b_i}\to R^{b_{i-1}}$
denotes the $(i-1)$st differential in this resolution. 

\begin{lemma}
\label{lemma1}
The sequence of homomorphisms: 
$$ 
C:\quad\cdots\to R^{b_3}\xrightarrow{d_3}R^{7}\xrightarrow{d_2}
R^{3}\xrightarrow{d_1}R^2\xrightarrow{d_0}R^2\xrightarrow{d_{-1}}R^2
\xrightarrow{d_{-2}} R^2\xrightarrow{d_{-3}} R^2\to\cdots 
$$ 
is a doubly infinite complex with $H_i(C)=0$ for all 
$i\in\mathbb Z$.
\end{lemma}

\begin{proof} (1). The defining equations of $R$ guarantee that
$d_{i-1}d_i=0$ for all $i\leq 2$. For $i\ge 2$ the maps $d_i$ are
differentials in a free resolution, hence the equality holds for all
$i\ge 3$, as well. We conclude that $C$ is a complex and $H_i(C)=0$
for all $i\ge 2$.

We let $(a,b)$ denote an element of $R^2$ written in the standard
basis of $R^2$ as a free $R$-module. For each $i\le 0$ the $k$-vector
space $\Ima d_{i}$ is generated by the elements:
\begin{align*}
d_i(1,0)&=(v,y)     &d_i(z,0)&=(vz,-vy)\\
d_i(0,1)&=(\alpha^{-i}x,z)  &d_i(0,v)&=(\alpha^{-i}vx,vz)\\  
d_i(v,0)&=(0,vy)  &d_i(0,x)&=(\alpha^{-i}vy,-\alpha^{-1}vx)\\
d_i(x,0)&=(vx,0) &d_i(0,y)&=(0,-vy)\\
d_i(y,0)&=(vy,-vx)    &d_i(0,z)&=(-\alpha^{-i-1}vx,0)
\end{align*}

It is clear that $\rank_k(\Ima d_{i})=8$ for all $i\le 0$. Since
$\rank_k(R^2)=16$, this implies that $\rank_k(\Ker d_i)=8$ for all
$i\le 0$, showing that $H_i(C)=0$ for all $i\leq -1$. 
 
Notice that for $i=1$ the elements above give $7$ linearly independent
elements in $\Ima d_1$, and the $8$th can be taken to be
$\varepsilon(0,0,1)=(yz,0)$. (Here $(a,b,c)$  denotes an element of 
$R^3$ in its standard basis as a free $R$-module.) Thus
$\rank_k(\Ima d_1)\ge 8$, and so $\rank_k(\Ker d_1)\leq 16$. In
particular, we obtain $H_0(C)=0$. 

To prove $H_1(C)=0$ we need to show that $\rank_k(\Ima d_2)\ge
16$. Indeed, the following elements in $\Ima d_2$ are linearly
independent:
\begin{align*}
d_2(e_1)&=(v,y,0)     &d_2(ye_4)&=(0,0,vy)\\
d_2(e_2)&=(\alpha^{-2}x, z,0) &d_2(ze_4)&=(0,0,vz)\\  
d_2(e_3)&=(-y,\alpha x,0) &d_2(xe_1)&=(vx,0,0 )\\
d_2(e_4)&=(0,0,v) &d_2(ye_1)&=(vy,-vx,0)\\
d_2(e_5)&=(0,0,x)    &d_2(ze_1)&=(vz, -vy,0)\\
d_2(e_6)&=(0,0,y) &d_2(ve_1)&=(0,vy,0)\\
d_2(e_7)&=(0,0,z)  &d_2(ve_2)&=(\alpha^{-2}vx, vz,0)\\
d_2(xe_4)&=(0,0,vx) &d_2(xe_2)&=(\alpha^{-2}vy,-\alpha^{-1}vx,0) 
\end{align*}
where $e_1,\dots ,e_7$ denote the elements comprising
the standard basis of $R^7$ as a free $R$-module. 
\end{proof}

If $f\colon M\to N$ is a homomorphism of $R$-modules, we let $f^*$
represent the induced map $\Hom_R(f,R)\colon \Hom_R(M,R)\to
\Hom_R(N,R)$. If $(D,\delta)$ is a complex of $R$-modules, 
then the complex $(D^*,\delta^*)$ has
$(D^*)_i=(D_{-i})^*$ and differentials $(\delta^*)_i=(\delta_{-i})^*$.
We write $\delta^*_i$ for $(\delta^*)_i$.

 Note that, upon identification of $R^*$ with $R$, the
map $d_i^*\colon R^2\to R^2$ for $i\geq 0$ is given in the standard basis 
of $R^2$ by the matrix
\[
\left(\begin{matrix}
v& y\\
\alpha^i x &z\end{matrix}\right).
\]
Similarly, the maps $d_{-1}^*$ and $d_{-2}^*$ are given by the
transposes of the matrices defining $d_1$ and $d_2$, respectively.

\begin{lemma} 
\label{lemma2}
The complex 
$$
C^*:\quad\cdots\to R^2\xrightarrow{d_2^*}
R^2\xrightarrow{d_1^*}R^2\xrightarrow{d_0^*}R^2
\xrightarrow{d_{-1}^*}R^3\xrightarrow{d_{-2}^*}R^7\to\cdots
$$ 
satisfies $H_i(C^*)=0$ if and only if $i\geq 1$.
\end{lemma}

\begin{proof} 
As a $k$-vector space, $\Ima d_i^*$ for $i\geq 0$ is generated by the
following elements:
\begin{align*}
d_i^*(1,0)&=(v,\alpha^i x)     &d_i^*(z,0)&=(vz,-\alpha^{i-1}vx)\\
d_i^*(0,1)&=(y,z)  &d_i^*(0,v)&=(vy,vz)\\  
d_i^*(v,0)&=(0,\alpha^ivx)  &d_i^*(0,x)&=(0,-\alpha^{-1}vx)\\
d_i^*(x,0)&=(vx,\alpha^ivy) &d_i^*(0,y)&=(-vx,-vy)\\
d_i^*(y,0)&=(vy,0)    &d_i^*(0,z)&=(-vy,0)
\end{align*}
One can see therefore that $\rank_k(\Ima d_i^*)=8$ if $i\geq 1$ and
$\rank_k(\Ima d_0^*)=7$. This implies that $\rank_k(\Ker d_i^*)=8$ if
$i\geq 1$ and $\rank_k(\Ker d_0^*)=9$, and it follows that $H_i(C^*)=0$ for
all $i\geq 1$ and $H_0(C^*)\ne 0$. 

For the proof that $H_i(C^*)\ne 0$ for $i=-1,-2$, note that the image of
$d^*_{-1}$ is generated as a $k$-vector space by the following elements.
\begin{align*}
d_{-1}^*(1,0)&=(v,\alpha^{-1}x,yz)     &d_{-1}^*(z,0)&=(vz,-\alpha^{-2}vx,0)\\
d_{-1}^*(0,1)&=(y,z,0)  &d_{-1}^*(0,v)&=(vy,vz,0)\\  
d_{-1}^*(v,0)&=(0,\alpha^{-1}vx,0)  &d_{-1}^*(0,x)&=(0,-\alpha^{-1}vx,0)\\
d_{-1}^*(x,0)&=(vx,\alpha^{-1}vy,0) &d_{-1}^*(0,y)&=(-vx,-vy,0)\\
d_{-1}^*(y,0)&=(vy,0,0)    &d_{-1}^*(0,z)&=(-vy,0,0)
\end{align*}
One sees easily that $\rank_k(\Ima d_{-1}^*)\le 8$.  
Therefore $\rank_k(\Ker d_{-1}^*)\ge 8$, and so $H_{-1}(C^*)\neq 0$.
Clearly $\rank_k(\Ker d_{-2}^*)$ consists of at least nine linearly independent
elements, namely the nine quadric elements in $R^3_2$. This shows that 
$H_{-2}(C^*)\neq 0$. 

Finally, we note from the matrix representing $d_2$ that 
$\Coker d_2\cong N\oplus k$, for some finite $R$-module $N$.  Therefore 
$H_i(C^*)\cong\Ext^{-i-2}_R(N\oplus k,R)\neq 0$ for all $i\leq -3$, since $R$ 
is not Gorenstein.
\end{proof}

\begin{chunk}
\label{moduleM} For each integer $s\ge 1$ let $M_s$ be the
cokernel of the map $d_{-s}\colon R^2\to R^2$. Using Lemma
\ref{lemma1}, note that
   $$M_s=\Coker(d_{-s})\cong\Ima(d_{-s-1})=\Ker(d_{-s-2})$$
and a truncation of the complex $C$ gives the beginning a minimal free
resolution of the $R$-module $M_s$: $$ \cdots\to R^7\xrightarrow{d_2}
R^3\xrightarrow{d_1}R^2\xrightarrow{d_0}R^2\xrightarrow{d_{-1}}R^2\to\cdots\to
R^2\xrightarrow{d_{-s+1}}R^2\xrightarrow{d_{-s}}R^2\to M_s\to 0$$ The
proof of Lemma \ref{lemma1} shows that $M_s$ has Hilbert series
$H_{M_s}(t)=2t+6t^2$. 
\end{chunk}

\medskip

We are now ready to state our main theorem:

\medskip

\begin{theorem}\label{main} For the family of $R$-modules
$\{M_s\}_{s\geq 1}$ defined above, we have:
\begin{enumerate}[\quad\rm(1)]
\item $M_s$ satisfies $({\mathsf{TR}_i})$ if and only if $i<s$. 
\item 
$\Tr(M_s)$ satisfies $(\mathsf{TR}_i)$ if and only if  $i>-s$.
\end{enumerate}
\end{theorem}

Note that this contains the Theorem stated in the introduction:
indeed, one can take there the modules $M_s$ to be as above and
$L=\Tr(M_1)$.

\begin{proof}[Proof of Theorem {\rm \ref{main}}] 
The second part of the theorem follows from the first part and the 
simple fact that a 
finite $R$-module $N$ satisfies the condition $(\mathsf{TR}_i)$ if and only if 
$\Tr(N)$ satisfies $(\mathsf{TR}_{-i})$.  
  
To compute $\Ext_R^*(N,R)$ for an $R$-module $N$ we take a
  minimal free resolution of $N$, we apply $(-)^*$ to it, and then
  compute homology of the resulting complex.

Applying $(-)^*$ to the minimal free resolution of $M_s$ given in
  \ref{moduleM}, and identifying $R$ with $R^*$, one obtains the complex
$$
\quad  R^2\xrightarrow{d_{s}^*}R^2\xrightarrow{d_{s-1}^*}R^2\to\cdots\to
R^2\xrightarrow{d_1^*}R^2\xrightarrow{d_{0}^*}R^2\xrightarrow{d_{-1}^*}R^3
\xrightarrow{d_{-2}^*}R^7\to\cdots
$$
Lemma \ref{lemma2} shows that $\Ext^i_R(M_s,R)=0$ for all $1\leq i\leq s-1$,
and $\Ext^i_R(M_s,R)\ne 0$ for $i\geq s$.  

A minimal free resolution of $\Tr(M_s)$ is given by 
$$
\cdots\to R^2\xrightarrow{d_{s+1}^*}R^2\xrightarrow{d_s^*}R^2,
$$
and applying $(-)^*$ we get
$$
R^2\xrightarrow{d_{-s}}R^2\xrightarrow{d_{-s-1}}R^2\to\cdots.
$$
Lemma \ref{lemma1} shows that $\Ext_R^i(\Tr(M_s),R)=0$ for all $i\ge 0$. 
This establishes (1), and hence the entire theorem.
\end{proof}

\section{Dependence}

Theorem \ref{main} shows that the conditions $(\mathsf{TR}_i)$ are, to
a large extent, independent.  However, in the Artinian graded case, there is some
overlap between these conditions.

Let $R$ be a standard graded ring and $M$ a finite graded $R$-module.
As noted by Avramov and Martsinkovsky in \cite{AM}, the module $M$ is
totally reflexive if and only if it satisfies $(\mathsf{TR}_i)$ for
all $i\ne -1$ if and only if it satisfies $(\mathsf{TR}_i)$ for all
$i\ne -2$. Thus, the condition $(\mathsf{TR}_i)$ for $i=-1$ is a 
consequence of the condition $(\mathsf{TR}_i)$ for all $i\neq-1$,
and the condition $(\mathsf{TR}_i)$ for $i= -2$ is a consequence 
of the condition $(\mathsf{TR}_i)$ for all $i\neq -2$.

The result above is based on a formula obtained by Avramov,
Buchweitz and Sally in \cite{ABS}. Buchweitz pointed out
to us that the same formula also yields the following.

\begin{proposition}\label{parity} Assume that $R$ is an Artinian 
standard graded
ring, and $M$ is a finitely generated graded $R$-module.  Let $A$ be a
finite set of integers of the same parity. If $M$ satisfies
$(\mathsf{TR}_i)$ for all $i\in \mathbb Z\smallsetminus A$, then the
module $M$ is totally reflexive.
\end{proposition}

\begin{proof} We may assume $0\notin A$. Suppose that $M$ satisfies
$(\mathsf{TR}_i)$ for all $i\in
\mathbb Z\smallsetminus A$. Since $A$ is a finite set, we conclude  that
$M$ satisfies $(\mathsf{TR}_i)$ for all $i$ with $|i|\gg 0$. This
allows us to use the main formula in \cite{ABS} which asserts equalities of
rational functions
\begin{align*}
\sum_{n\in \{0\}\cup A}(-1)^nH_{\Ext^n_R(M,R)}(t)&=
\frac{H_M(t^{-1})\cdot H_R(t)}{H_R(t^{-1})}\quad   &(1)\\
\sum_{-n\in\{0\}\cup A}(-1)^nH_{\Ext^n_R(M^*,R)}(t)&=
\frac{H_{M^*}(t^{-1})\cdot H_R(t)}{H_R(t^{-1})}\quad  &(2) 
\end{align*}

We set 
$$ 
P(t)=\sum_{n\in A}H_{\Ext^n_R(M,R)}(t)\quad \text{ and
}\quad Q(t)=\sum_{-n\in A}H_{\Ext^n_R(M^*,R)}(t). 
$$ 

Let $\sigma=0$ if $A$ contains only odd integers and 
$\sigma=1$ if $A$ contains only even integers.
The formulas (1) and (2) give then
\begin{gather*}
H_{M^*}(t)=\frac{H_M(t^{-1})\cdot H_R(t)}{H_R(t^{-1})}+(-1)^{\sigma}P(t)\\H_{M^{**}}(t)=\frac{H_{M^*}(t^{-1})\cdot
H_R(t)}{H_R(t^{-1})}+(-1)^{\sigma}Q(t)
\end{gather*}
Substituting the formula for $H_{M^*}(t^{-1})$ given by the first
equation into the second equation, we obtain: $$
H_{M^{**}}(t)=H_M(t)+\frac{ H_R(t)}{H_R(t^{-1})}\cdot
(-1)^{\sigma}P(t^{-1})+(-1)^\sigma Q(t) $$ and it follows that $$
H_R(t)\cdot P(t^{-1})+H_R(t^{-1})\cdot
Q(t)=(-1)^{\sigma}H_R(t^{-1})\big(H_{M^{**}}(t)-H_M(t)\big).  $$

Assume that $A$ contains only odd integers.  We have then
$-2\notin A$, hence $(\mathsf{TR}_i)$ is satisfied for $i=-2$. The map
$M\to M^{**}$ is then surjective, implying a coefficientwise
inequality $H_{M^{**}}(t)\le H_M(t)$. Since $\sigma=0$ in this case
and $H_R(t^{-1})$ has positive coefficients, it follows that the
Laurent polynomial on the right has nonpositive coefficients.

  Both terms of the left-hand side sum are Laurent polynomials with
  nonnegative coefficients, and it follows that $P(t)=0$ and $Q(t)=0$,
  implying that $\Ext^n_R(M,R)=0$ for all $n\in A$ and
  $\Ext^n_R(M^*,R)=0$ for all $n$ with $-n\in A$. In conclusion,
  $(\mathsf{TR}_i)$ is satisfied for all $i\ne -1$.
Furthermore, we conclude from the formula above that
$H_{M^{**}}(t)=H_M(t)$.  Since the map
  $M\to M^{**}$ is surjective, it follows that it is an isomorphism,
  hence $(\mathsf{TR}_i)$ is satisfied for $i=-1$ as well.

Proceed similarly when $A$ contains only even integers. 
\end{proof}

Proposition \ref{parity} leads to the following question:

\begin{question}
Let $R$ be a commutative (local) Artinian ring. If a finite
$R$-module $M$ satisfies $(\mathsf{TR}_i)$ for all but finitely many
values of $i$, does it follow that $M$ is totally reflexive? 
\end{question}

\section{Minimimal acyclic complexes of free modules}

Let $S$ be a commutative Noetherian local ring with
maximal ideal $\fn$.

A complex $F$ of free $S$-modules
$$
\cdots\to F_{i+1}\xrightarrow{\phi_{i+1}} F_i\xrightarrow{\phi_i} F_{i-1}\to \cdots\ 
$$ is said to be {\it minimal} if $\phi_i(F_i)\subseteq \fn F_{i-1}$
for all $i$. The complex $F$ is \emph{acyclic} if $H_i(F)=0$ for all
$i\in\mathbb Z$. For any minimal acyclic complex of finite free
$S$-modules $F$ we can consider two sequences:
$${\bd
\beta}_F^+:=\{\rank_RF_i\}_{i\ges 0}\quad\text{ and}\quad
{\bd\beta}^-_F:=\{\rank_RF_{-i}\}_{i\ges 0}$$ Assuming that
$F_i\ne 0$ for all $i$, it is then natural to ask whether these
two sequences have similar asymptotic behavior.

A sequence $\{\beta_i\}_{i\ge 0}$ is said to have {\it exponential
growth} if there exist numbers $1<A\le B$ such that inequalities
$A^i\le \beta_i\le B^i$ hold for all $i\gg 0$.

When the maximal ideal of $S$ satisfies $\fn^3=0$, Lescot \cite{L}
proved that the Betti numbers of a finitely generated
$S$-module $N$ are either eventually stationary, or they have
exponential growth; in the last case they are eventually strictly
increasing. It is clear from the Poincar\'e series given in
\ref{ringR} that the Betti numbers of $k$ over our ring $R$ have
exponential growth. Furthermore, with $d_2$ as defined there, since
$\Coker(d_2)$ has a copy of $k$ as a direct summand, its Betti numbers
have exponential growth and are eventually strictly increasing. 

\medskip

In conclusion, the complex $C$ of Lemma \ref{lemma1} has the following
properties:
\begin{enumerate}[\quad (a)]
\item  ${\bd \beta}_C^+$ has exponential
  growth and is eventually strictly increasing.
\item   ${\bd\beta}^-_C$ is constant (nonzero).
\end{enumerate}

\smallskip

Several questions arise:
\begin{question}
  Does there exist a ring $S$ as above and a minimal acyclic complex of
  free nonzero $S$-modules $F$ such that ${\bd \beta}_F^-$ has exponential
  growth (or is eventually strictly increasing) and ${\bd\beta}^+_F$ is
  eventually constant?
\end{question}

\begin{question}
  Do there exist examples of different asymptotic behavior for ${\bd
  \beta}^-_F$ and ${\bd\beta}^+_F$ if we also require $H(F^*)=0$ ?
  Can such examples exist over a Gorenstein ring?
\end{question}

The last question is equivalent to asking whether the Betti numbers of $M$ and
$M^*$ can have different asymptotic behavior when $M$ is totally
reflexive, and, in particular, when $S$ is Gorenstein. Theorem 5.6 of
\cite{AB} shows that the answer to this question is ``no'' when $S$ is a
complete intersection.

\section*{Acknowledgements}

The authors thank Luchezar Avramov, Ragnar Buchweitz, Lars
Christensen, Srikanth Iyengar, and Hamid Rahmiti 
for their useful comments and suggestions.

\end{document}